\documentstyle{amsppt}
\define\2x2#1#2#3#4{\pmatrix #1&#2\\#3&#4\endpmatrix}
\define\a#1#2{A_{#1#2}}
\define\be#1#2{B_{#1#2}}
\define\ce#1#2{C_{#1#2}}
\define\amod{\g A\text{\bf mod}}

\define\amodb{\g A\text{\bf mod}\g B}
\define\amoda{\g A\text{\bf mod}\g A}

\define\moda{\text{\bf mod}\g A}

\define\bmoda{\g B\text{\bf mod}\g A}
\define\M#1#2{M^{#1#2}}
\define\m#1#2{m^{#1#2}}

\define\g{\goth}
\define\A#1{\Cal A(#1)}
\define\F#1{\Cal F(#1)}
\define\N#1{\Cal N(#1)}
\define\I#1{\Cal I(#1)}
\define\Si#1{\Cal S(#1)}
\define\T#1{#1\widehat\otimes#1^*}
\define\Bn#1{\Cal B(#1)}
\define\Aa#1#2{\Cal A(#1,#2)}
\define\Ff#1#2{\Cal F(#1,#2)}
\define\nul{\{0\}}
\define\flp#1{\flushpar{\rm(#1)}}
\define\Ii#1#2{\Cal I(#1,#2)}
\define\Nn#1#2{\Cal N(#1,#2)}
\define\Bb#1#2{\Cal B(#1,#2)}
\define\ot#1{\underset#1\to{\otimes}}
\define\pot#1{\underset#1\to{\widehat\otimes}}
\define\pott{\widehat\otimes}
\define\pota{\pot{\g A}}
\define\potb{\pot{\g B}}
\define\pra#1#2{\{#1,#2\}}
\define\prb#1#2{[#1,#2]}
\define\tens#1{#1\pot{#1}#1}
\define\btens{\bar\otimes}
\define\Ho#1#2{\Cal H^{#1}(#2,#2^*)}
\define\HoM#1{\Cal H^1_{\g M}(#1,#1^*)}
\define\apb{{}_{\g A}P_{\g B}}
\define\bqa{{}_{\g B}Q_{\g A}}

\define\lp#1{\ell_p(#1)}
\define\LP#1{L_p(\mu,#1)}
\define\nrm#1{\Vert #1\Vert}

\define\si{self-induced}
\define\tr{\operatorname{tr}}

\define\cl{\operatorname{cl}}
\define\clspan{\operatorname{clspan}}
\define\way{weak amenability}
\define\wae{weakly amenable}
\redefine\l{\langle}
\define\r{\rangle}
\define\set#1{\{#1\}}
\refstyle{A} 
\topmatter
\title Self-induced Banach algebras\endtitle
\author Niels Gr\o nb\ae k \endauthor
\address Department of Mathematics, Institute for Mathematical
Sciences, Universitets\-parken 5, DK-2100
Copenhagen {\O}, Denmark\endaddress
\email gronbaek\@math.ku.dk\endemail
\abstract A Banach algebra $\g A$ is \si\ if the
multiplication $\g A \pot{\g A}\g A\mapsto\g A$ is an isomorphism. The
class of \si\ Banach algebras is a natural generalization of
unital Banach algebras, providing a fertile framework for developing
homological aspects of unital Banach algebras. Elementary results with
applictions to computations of the bounded cohomology groups $\Ho1{\g A}$,
with emphasis on $\g A=\A X$, the approximable operators on a Banach
space $X$, are given.
\endabstract
\endtopmatter
\document

\subheading{0. Introduction} A Banach algebra $\g A$ is {\it \si} if
the multiplication $$\g A\pot{\g A}\g A@>\mu_{\g A}>>\g A$$ is an
isomorphism. The concept of \si\  Banach algebras was introduced in [G2]
with the purpose of providing a setting for a Banach algebra theory of
Morita equivalence rich enough to include commonly occurring Banach
algebras. Self-inducedness is in liking with {\it H-unitality}, the
concept introduced by M\. Wodzicki in [W]. The Banach algebraic
version of H-unitality is stipulated in the requirement that the bar
complex of $\g A$ 
$$
0@<<<\g A@<<<\g A\pot{}\g A@<<<\g A\pot{}\g A\pot{}\g A@<<<\cdots
$$
is pure-exact, or equivalently, the bounded homology
groups $\Cal H_n(\g A,\g X)=\{0\},\;n\geq1$ for any annihilator Banach $\g
A$-bimodule $\g X$. It is rather elementary that if $\g A$ is H-unital as
a Banach algebra, then the inclusion  $\g A\mapsto\g A_+$,
where $\g A_+$ is the unitization of $\g A$, implements an isomorphism 
$\Ho n{\g A_+}\cong\Ho n{\g A}$. Rephrasing this, the admissable short
exact sequence $0@>>>\g A@>>>\g A_+@>>>\Bbb C@>>>0$ implements a long
exact sequence
$\cdots@>>>\Ho n{\Bbb C}@>>>\Ho n{\g A_+}@>>>\Ho n{\g
A}@>>>\cdots$. However, the main result of [W]  that this excission
property holds in general, is far from elementary: If $\g A$ is
H-unital as a Banach algebra, then any weakly admissable extension
$0@>>>\g A@>>>\g B@>>>\g C@>>>0$ implements a long exact sequence
$$
\cdots@>>>\Ho n{\g C}@>>>\Ho n{\g B }@>>>\Ho n{\g A }@>>>\cdots\;.
$$

Being \si\  may be seen as H-unitality in ``degree $1+\frac12$''. More precisely
$\g A$ is \si\  is equivalent to each of 
\roster
\item $\Cal H^1(\g A,\g X)=\Cal H^2(\g A,\g X)=\{0\}$ for any annihilator
bi-module $\g X$.
\item $\Cal H^0_{\text{bar}}(\g A)=\Cal H^1_{\text{bar}}(\g A)=\{0\}$,
where $\Cal H^{\bullet}_{\text{bar}}$ is the cohomology of the dual bar-complex.
\endroster

Besides providing the framework for a Morita theory, \si\ Banach
algebras share a basic cohomological property with unital algebras: If
$\g A@>>>\g B$ is an embedding of a \si\ Banach algebra $\g A$ as an
ideal in $\g B$, then  restriction implements isomorphisms $\Cal
H^1(\g A,\g B^*)\cong H^1(\g B,\g A^*)\cong H^1(\g A,\g A^*)$, see
[G3, Lemma 2.2]. In the present paper we show some elementary properties of \si\
Banach algebras with the aim of computing the cohomology groups
$\Ho1{\g A}$. An imporpant aspect will be Morita invariance of Hochschild
cohomology for \si\ Banach algebras. We shall apply our findings to
algebras of the type $\A X$, the approximable operators on a Banach
space $X$, thereby improving on previous results, in particular
concerning direct sums and tensor products, of [B1], [B2], and [G3].

\subheading{1. Preliminaries} For  Banach spaces $X$ and $Y$ we
consider the following spaces of operators 
$$
\align
&\Ff XY=\{\text{finite rank operators  $X\to Y$}\}\\
&\Aa XY=\{\text{approximable operators  $X\to Y$}\}=\overline{\Ff XY}\\ 
&\Ii XY=\{\text{integral operators  $X\to Y$}\}\\
&\Bb XY=\{\text{bounded operators  $X\to Y$}\} 
\endalign
$$
As customary we shall write $\operatorname{Operators}(X)$ for
$\operatorname{Operators}(X,X)$. 

We shall use $\l\cdot\,,\cdot\r$ to denote the bilinear form $\l
x,x^*\r=x^*(x)$ for $x\in X,x^*\in X^*$. The corresponding bilinear
functional $\tr\colon X\otimes X^*\to\Bbb C$ is {\it the canonical
trace}. Since it {\it is} canonical it is not annotated with the space
$X$.

The definitions of  Banach (co)homological concepts are standard and
can be found for example in [H] and [J].
For Banach algbras $\g A$ and
$\g B$ we use the notations $\amod$, $\moda$, and $\amodb$ for the
categories of Banach left, Banach right $\g A$-modules, and Banach  $\g
A$-$\g B$-bimodules.  The unitization of $\g A$ is the Banach algebra
given by the augmentation $\g A@>>>\g A_+@>>>\Bbb C$. We denote bounded Hochschild
homology and cohomology for $\g A$ with coefficients in $\g X\in
\amoda$ by $\Cal H_n(\g A,\g X)$ and $\Cal H^n(\g A,\g X)$. We shall
almost exclusively be interested in $n=1$.  

Recall

\definition{1.1 Definition} Let $\g A$ be a Banach algebra, let
$P\in\moda$  and  $Q\in\amod$. We define
$$
P\pot{\g A}Q=P\pot{}Q/N
$$
where $\pot{}$ is the projective tensor product and
$N=\clspan\{p.a\otimes q-p\otimes a.q\mid p\in P,q\in Q,a\in \g 
A\}$. Thus, $P\pot{\g A}Q$ is the universal object for linearizing
bounded, $\g A$-balanced bilinear maps $P\times Q\to Z$ into Banach
spaces $Z$.
\enddefinition

For $P\in\amod$ the we denote the linear map implemented by multiplication
by $\mu_P\colon\g A\pot{\g A}P@>>>P$ and likewise for right and
two-sided modules. If it is unambiguous we
shall omit the subscript.

All (multi)-linear maps will be assumed to be bounded, and all
occuring series in Banach spaces are assumed to be absolutely convergent.
  
\subheading{2. Morita contexts and derivations of \si\ Banach algebras}

\definition{2.1 Definition} A Banach algebra $\g A$ is called {\it
  \si } if multiplication implements a bimodule isomorphism 
$$
\g A\pot{\g A}\g A@>\mu>>\g A.
$$
More generally, $P\in\amod$ is an {\it $\g A$-induced module} if $\mu_P$
is an isomorphism, and similarily for $P\in\moda$.
\enddefinition

As mentioned, one motivation for this definition is that if $\g A$ is
\si, then $\Ho1{\g A_+}\cong\Ho1{\g A}$. However, \si ness is not
necessary for this. More precisely we have

\proclaim{2.2 Proposition} The following three statements are equivalent.

\flp i $\Ho1{\g A_+}\cong\Ho1{\g A}$;
\flp{ii} $\Cal H^1(\g A,\Bbb C)=\{0\}$ and $\Cal H^2(\g A,\Bbb
C)\to\Cal H^2(\g A,\g A_+^*)$ is injective for the annihilator module
$\Bbb C$;
\flp{iii} $\cl(\g A^2)=\g A$ and for each derivation $D\colon\g
A\to\g A^*$ the bilinear map $(a,b)\mapsto\l a,D(b)\r+\l b,D(a)\r$ is
$\g A$-balanced.
\endproclaim
\demo{Proof} This is proved by inspection of the concepts involved.
\enddemo

\example{2.3 Examples} (a) $\g A$ is \si, if $\g A$ is left or right
flat and $\cl(\g A^2)=\g A$, in particular if $\g A$ is
biflat. Important instances occur, when $\g A$ has a one-sided bounded
approximate identity, in which case $\g A$ is in fact H-unital.
\flp{b} For a Banach space $X$ the Banach algebra $\A X$ is \si,
\roster
\item if $X$ has the bounded approximation property, cf\. (a), or more
generally if $\mu_{\A X}$ is surjective and $X$ has the approximation
property, see [G3];
\item if $X=\ell_2(E)$, when $E$ and $E^*$ both are of cotype 2,
[Proposition 4.2, B1];
\item if $X=E\oplus C_p,\;1\leq p\leq\infty$ with $E$ arbitrary and
$C_p$ one of the spaces defined by W.B\. Johnson, [Jo].
\endroster
\flp c If $P$ is one of the spaces constructed by G\. Pisier in [P1],
then $\mu_{\A P}$ is surjective, but not injective.
\flp d If $P$ is the space constructed by G\. Pisier in [P2], then
$\mu_{\A P}$ is not surjective.
\flp{e} $\N X$ is \si, if and only if $X$ has the approximation
property, [G3].
\flp{f} The augmentation ideal of the free
group on two generators, $I(\Bbb F_2)$, is not \si, yet
$\Ho1{I(\Bbb F_2)}=\nul$, [G\&L, Corollary 3.2].
\endexample

In the setting of this paper it will be convenient to have Morita
equivalence described in terms of Morita contexts. For Banach algebras
with trivial annihilator this is equivalent
to the definition given in [G2].

\definition{2.4 Definition} Let $\g A$ and $\g B$ be \si\ Banach algebras,
and let $P\in\amodb$ and $Q\in\bmoda$. Two balanced bilinear maps
$\pra{~}{}\colon P\times Q\to\g A$ and $\prb~~\colon
Q\times P\to\g B$ are called {\it compatible pairings} if they implement
bimodule homomorphisms $P\potb Q\to\g A$ and $Q\pota P\to\g B$ and
$$
\pra pqp'=p\prb q{p'},\;\prb qpq'=q\pra p{q'};\quad p,p'\in
P,\;q,q'\in Q.\tag\S 
$$ 

These data are collected in  a {\it Morita context\/} $\g M(\g A,\g
B,P,Q,\pra~{} ,\prb~~)$. To $\g M$ is associated a Banach algebra $\Cal
M$ consisting of $2\times2$ matrices
$$
\2x2 apqb;\quad a\in \g A,b\in \g B,p\in P,q\in Q
$$
with the product defined by means of module multiplication and the
compatible pairings. The Morita context is called {\it full\/} if both
pairings implement bimodule isomorphisms. The Banach algebras $\g A$ and
$\g B$ are {\it Morita equivalent}, if there are modules
$P\in\amodb,\;Q\in\bmoda$ and compatible pairings $\pra~{},\;\prb~~$
constituting a full Morita context.
\enddefinition

The simplest case of Morita invariance of Hochshild cohomology occurs
for $n\times n$-matrices. First we investigate \si ness:

\proclaim{2.5 Proposition} Let $\g A$ be a Banach algebra satisfying
$\cl(\g A^2)=\g A$. Then 
$$
M_n(\g A)\pot{M_n(\g A)}M_n(\g A)=M_n(\tens{\g A})\;,
$$
where $M_n$ denotes $n\times n$ matrices.
\endproclaim
\demo{Proof} We prove that $M_n(\tens{\g A})$ has the universal
property defining the tensor product $\pot{M_n(\g A)}$. Let
$$
\Phi\colon M_n(\g A)\times M_n(\g A)\to Z
$$
be a bounded $M_n(\g A)$-balanced bilinear map into a Banach space $Z$. Let
$E_{ij}$ be the elementary matrices in $M_n(\Bbb C)$ and let
$a,a',b\in\g A$. Then
$$
\align
\Phi(aa'E_{ij},bE_{kl})&=\Phi(aE_{i1}a'E_{1j},bE_{kl})=\Phi(aE_{i1},a'E_{1j}bE_{kl})\\
&=\cases
0,&j\neq k\\
\Phi(aE_{i1},a'bE_{1l}),&j=l
\endcases
\endalign
$$
Since $\cl(\g A^2)=\g A$, we further get
$\Phi(aE_{ij},bE_{kl})=0$ for $j\neq k$. Let 
$$
\phi_{ij}(a,b)=\Phi(aE_{i1},bE_{1j})
$$
Then $\phi_{ij}$ is clearly $\g A$-balanced. Let 
$\iota\colon M_n(\g A)\times M_n(\g A)\to M_n(\tens{\g A})$ be
the bilinear map
$$
((a_{ij}),(b_{kl}))\mapsto(\sum_ka_{ik}\ot{\g A}b_{kj})_{ij}.
$$
Then $\iota$ is  $M_n(\g A)$ balanced and the above
calculation shows that if we for $U=(u_{ij})\in M_n(\tens{\g A})$ set
$$
\widetilde\Phi(U)=\sum_{ij}\phi_{ij}(u_{ij})
$$
then $\Phi=\widetilde\Phi\circ\iota$. Since $\operatorname{span}\iota(M_n(\Cal
A)\times M_n(\g A))$ is dense, $\widetilde\Phi$ is uniquely
determined.
\enddemo
\medskip

\proclaim{2.6 Corollary} $\g A$ is \si\ if and only if $M_n(\Cal
A)$ is \si .
\endproclaim

\smallskip
We now have elementary Morita invariance of Hochschild
cohomology in degree 1.

\proclaim{2.7 Theorem} Let $\g A$ be \si\, let
$\iota\colon\g A\pot{}\g A\to M_n(\g A)\pot{}M_n(\g A)$
be the map given by
$$
a\mapsto\pmatrix a&\\&{\matrix 0&&\\&\ddots&\\&&0\endmatrix} \endpmatrix,
$$
and let $\tau\colon M_n(\g A)\pot{}M_n(\g A)\to \g A\pot{}\g A$ be the
map
$$
\tau((a_{ij})\otimes(b_{ij}))=\sum_{i,j}a_{ij}\otimes b_{j\,i}.
$$
Then the implemented maps $\Ho1{M_n(\g A)}@>\iota^*>>\Ho1{\g A}$ and
$\Ho1{\g A}@>\tau^*>>\Ho1{M_n(\g A)}$ are inverses of one another.
\endproclaim
\demo{Proof} If $\g A$ is unital, this is the usual Morita invariance of
Hochschild cohomology of matrices. Consider the admissable sequence
$$
\nul\to M_n(\g A)\to M_n(\g A_+)\to M_n(\Bbb C)\to\nul
$$
We have $\Cal H^1(\g A_+,\Bbb C)=\Cal H^2(\g A_+,\Bbb C)=\nul$, since
$\g A$ is \si. By Morita invariance we then have $\Cal H^1(M_n(\g A_+),M_n(\Bbb
C))=\Cal H^2(M_n(\g A_+),M_n(\Bbb C))=\nul$. Since $M_n(\g A)$ is \si,
we also have $\Cal H^1((M_n(\g A),M_n(\Bbb
C))=\Cal H^2(M_n(\g A),M_n(\Bbb C))=\nul$ for the annihilator module
$M_n(\Bbb C)$, cf\. the introduction.

>From the long exact sequence of cohomology we get a commutative diagram
$$
\CD
\nul@>>> \Ho1{M_n(\g A_+)}@>>>\Cal H^1(M_n(\g A_+),M_n(\g A)^*)@>>>\nul\\
&&@VVV@VVV&\\
\nul@>>> \Cal H^1(M_n(\g A),M_n(\g A_+)^*)@>>>\Cal H^1(M_n(\g A),M_n(\g A)^*)@>>>\nul
\endCD
$$
where the vertical maps are implemented by inclusions. The rows are exact
and the right vertical arrow is an isomorphism by the extension property of
\si\ Banach algebras ([G3,Lemma 2.2]). Hence the left vertical arrow is an
isomorphism, and we have a commutative diagram
$$
\CD
\Ho1{M_n(\g A_+)}@>>\cong>\Ho1{M_n(\g A)}\\
@AA\iota^*A@AA\iota^*A\\
\Ho1{\g A_+}@>>\cong>\Ho1{\g A}
\endCD
$$
It follows than $\iota^*\colon\Ho1{M_n(\g A)}\to\Ho1{\g A}$ is an
isomorphism.
\enddemo

We want to pursue this further for more general Morita contexts. First
we need to  precisionize right exactness of the tensor functor:

\proclaim{2.8 Lemma} Consider a short complex in $\amod$
$$
X@>i >>Y@>q>>Z\to\nul
$$
with $\cl(i(X))=\ker(q)$ and $q(Y)=Z$ (i.e\. the dual complex
$X^*@<i^* <<Y^*@<q^*<<Z^*@<<<\nul$ is exact). Then for any $P\in\moda$:
$$
\cl((\bold1_P\otimes i)(P\pot{\g A}X))=\ker(\bold1_P\otimes q)\text{
and }(\bold1_P\otimes q)(P\pot{\g A}Y)=(P\pot{\g A}Z).
$$
\endproclaim
\demo{Proof} The proof of [G1, Lemma 3.1] works verbatimly.
\enddemo

We use this to weaken the defining properties of Morita equivalence:

\proclaim{2.9 Lemma} Let $\g A$ and $\g B$ be any Banach
algebras, and let $P\in\amodb$ and $Q\in\bmoda$. Suppose that 
\flp1 $\g A$ is \si;
\flp2 there are pairings
\roster
\item"" $\pra\ \ \colon P\times Q\to\g A$ 
\item"" $\prb\ \ \colon Q\times P\to\g B$ 
\endroster
\phantom{\flp2} satisfying the compatibility conditions {\rm(\S)};
\flp3 $\pra\ \ $ implements an epimorphism;
\flp4 $\cl(\g A.P)=P$.

\smallskip\noindent Then $\pra\ \ $ actually implements an isomorphism.
\endproclaim
\demo{Proof} Let $q\colon P\pot{\g B}Q\to\g A$ be the epimorphism implemented by
$\pra\ \ $ and set $K=\ker q$, so that we have a short exact sequence
$$
\nul\to K@>i>>P\pot{\g B}Q@>q>>\g A\to\nul.
$$
We want to show that $K=\nul$. Consider 
$$
\CD
\g A\pot{\g A}K@>\bold1_{\g A}\ot{\g A}i>>\g A\pot{\g A}P\pot{\g
B}Q@>\bold1_{\g A}\ot{\g A}q>>\g A\pot{\g A}\g A@>>>\nul\\
@V\mu_1VV@V\mu_2VV@V\mu_3VV\\
K@>i>>P\pot{\g B}Q@>q>>\g A@>>>\nul
\endCD
$$
where we for brevity have set $\mu_1=\mu_K,\; \mu_2=\mu_{P\pot{\g B}Q}$ and
$\mu_3=\mu_{\g A}$. By assumption $\mu_3$ is an isomorphism and
$\mu_2$ has dense range. We show that $\mu_1$ has dense range. Let $k\in
K$. Since $\mu_2$ has dense range, there is a sequence $t_n\in\g A\pot{\g A}P\pot{\g
B}Q$ so that $\lim \mu_2(t_n)=i(k)$. Hence $q(\mu_2(t_n))\to
q(i(k))=0$, so $\mu_3^{-1}\circ q\circ\mu_2(t_n)=\bold1_{\g
A}\ot{\g A}q(t_n)\to0$. By Lemma 2.8 and the open mapping theorem there is a
sequence $u_n\in\g A\pot{\g A}K$ so that $t_n-\bold1_{\g
A}\ot{\g A}i(u_n)\to0$. But then $\lim\mu_2(t_n-\bold1_{\g
A}\ot{\g A}i(u_n))\to0$, i.e\. $\lim\mu_2(\bold1_{\g
A}\ot{\g A}i(u_n))=i(k)\to0$ or equivalentlly $\mu_1(u_n)\to k$.

Let $k=\sum p_n\ot{\g B}q_n$. For any $p\in P,\; q\in Q$ we have
$$
\multline
\pra pq k=\sum \pra pq p_n\ot{\g B}q_n= \sum p\prb q{p_n}\ot{\g
B}q_n\\
=\sum p\ot{\g B}\prb q{p_n}q_n=\sum p\ot{\g B} q\pra{p_n}{q_n}=0
\endmultline
$$
so by (3) $\g A.K=\nul$ and we are done.
\enddemo

\proclaim{2.10 Theorem} Let $\g M(\g A,\g 
B,P,Q,\dots)$ be a Morita context of \si\ Banach algebras $\g A$ and $\g B$.
Suppose that the pairings both implement epimorphisms.
\flp1 With $\apb=\g A\pot{\g A}P\pot{\g B}\g B$ and $\bqa=\g B\pot{\g
B}Q\pot{\g A}\g A$ we have a full Morita context $\g M(\g A,\g
B,\apb,\bqa,\dots)$. In particular $\g A\sim\g B$
\flp2 Let $\Cal M$ be the Banach algebra associated with $\g M(\g A,\g
B,\apb,\bqa,\dots)$. Then $\Cal M$ is \si\ and $\Cal M\sim\g A\sim\g B$.
\endproclaim
\demo{Proof} We have a commutative diagram 
$$
\CD
P\pot{\g B}Q\pot{\g A}P\pot{\g B}Q@>1>>\tens{\g A}\\
@V2VV@V3VV\\
P\pot{\g B}\g B\pot{\g B}Q@>4>> \g A
\endCD
$$
Since 1, 2, and 3 are epimorphism, 4 is also an epimorphism, so we
have epimorphisms
$$
\align
&(\g A\pot{\g A}P\pot{\g B}\g B) \pot{\g B} (\g B\pot{\g B}Q\pot{\g
A}\g A) \to\g A\\
&(\g B\pot{\g B}Q\pot{\g A}\g A) \pot{\g A} (\g A\pot{\g A}P\pot{\g
B}\g B) \to\g B
\endalign
$$
It follows that $\apb$ and $\bqa$ satisfy the hypotheses of Lemma, thereby
proving (1).

Thus we may assume that we have a full Morita context $\g M(\g A,\g
B,P,Q,\pra~{},\prb~~)$ with associated Banach algebra
$$
\Cal M=\2x2{\g A}PQ{\g B}\,,
$$
by replacing the bimodules $P$ and $Q$ by the induced
bimodules $\apb$ and $\bqa$
if necessary. Define $\g P=(\g A\;P)$ and $\g Q=\pmatrix \g
A\\Q\endpmatrix$. Then 
$\g P$ and $\g Q$ are both $\g A$- and $\Cal M$- induced, and by means of matrix
multiplication we have epimorphic compatible pairings. In order to
finish by appealing to Lemma 2.9 we must show that $\Cal M$ is \si. It is
convenient to write $\Cal M=\2x2 {\g A_{11}}{\g A_{12}}{\g A_{21}}{\g
A_{22}}$. Let $\phi\colon\Cal M\times\Cal M\to \Bbb C$ be an $\Cal M$-
balanced bilinear functional. Using the isomorphisms $\g A_{ij}\pot{\g
A_{jj}}\g A_{jk}\cong \g A_{ik},\;i,j,k=1,2$ one finds linear functionals
$\varphi_{ij}\colon\g A_{ij}\to \Bbb C$ such that with
$\tilde\phi=\2x2{\varphi_{11}}{\varphi_{12}}{\varphi_{21}}{\varphi_{22}}$
one has $\phi((\g a_{ij}),(\g \alpha_{ij}))=\tilde\phi((\g a_{ij})(\g
\alpha_{ij}))$. This is a routine but rather lengthy exercise along
the lines of the proof of Proposition 2.5 and is therefore omitted. It follows
that $\Cal M$ has the universal property of $\tens{\Cal M}$.
\enddemo

\remark{Remark}
We cannot conclude the reverse to (2). Let $X$ be a Banach space such that
$\A X$ is not \si . Since $\A{X\oplus C_p}$ is \wae\ ([B1, Corollary 3.3] and
$X\oplus C_p$ is decomposable, then
$\A{X\oplus C_p}$ is \si , ([G3, Theorem 2.9].
\endremark
\medskip
In order to compare cohomology of Morita equivalent Banach algebras we
shall as previously exploit the double complex of Waldhausen [DI]. Let
$\g M(\g A,\g B,P,Q,\dots)$ be a full Morita context of \si\ Banach
algebras.  The complex on the first axis of the Waldhausen bi-complex
is the shifted Hochschild complex $C_{\bullet-1}(\g A,\g A)$ and on
the second axis the shifted Hochschild complex $C_{\bullet-1}(\g B,\g
B)$ The $n$'th column for $n\geq1$ is the complex $P\pot{\g
B}C_{\bullet}(\g B,Q\pot{}\Cal A^{\pot{}(n-1)})$ where $C_{\bullet}(\g
B,Q\pot{}\g A^{\pot{}(n-1)})$ is the normalized bar resolution of the
left $\g B$-module $Q\pot{}\Cal A^{\pot{}(n-1)}$. Similarly, the
$m$'th row is the complex $Q\pot{\Cal A}C_{\bullet}(\g A,P\pot{}\g
B^{\pot{}(m-1)})$. For details, see [G2].

The cohomology groups $\Ho n{(\cdot)}$ are related to the dual
bi-cocomplex:

$$
\CD
 (3)\;@. (\g B\pot{}\g B\pot{}\g B )^*@>>>\\
@.    @AAA     @AAA\\
 (2)@. (\g B\pot{}\g B)^*@>>>( P\pot{}\g B\pot{}Q )^*@>>>\\
@.    @AAA        @AAA   @AAA \\
 (1)@. \g B^*@>>>( P\pot{} Q )^*@>>>(P\pot{} Q \pot{}\g A )^*@>>>\\
@.  @AAA        @AAA @AAA        @AAA\\
@. 0@>>>\g A^*@>>>(\g A\pot{}\g A)^*@>>>(\g A\pot{}\g A\pot{}\g A)^*\\
@.@.  (1)@. (2)@.     (3)@.
\endCD
$$

Denote the $E_p-$terms for the spectral sequence whose $E_2$-terms are
obtained by taking
first horizontal and then vertical cohomology by $^{I}\!E^{(ij)}_p$ and 
the $E_p-$terms for the spectral sequence whose $E_2$-terms are obtained by taking
first vertical and then horizontal cohomology by $^{II}\!E^{(ij)}_p$. We
compute
$$
\matrix
^{I}\!E^{(i0)}_1=\Ho{i-1}{\g A} &^{I}\!E^{(ij)}_1=\nul,\quad i=0,1;\,j\geq 1\\
^{II}\!E^{(0j)}_1=\Ho{j-1}{\g B} &^{II}\!E^{(ij)}_1=\nul,\quad j=0,1;\,i\geq 1
\endmatrix
$$
where the terms on the axes follow from the definition of Hochschild
cohomology and the others from the Morita context being full and the
definition of the tensor products $\pot{\g A}$ and $\pot{\g B}$. 

Thus
$^{I}\!E^{(20)}_p$ and $^{II}\!E^{(02)}_p$ stabilize at $p=2$ and
$^{I}\!E^{(02)}_1=^{I}\!E^{(11)}_1={}^{I}\!E^{(20)}_1={}^{II}\!E^{(20)}_1=
{}^{II}\!E^{(11)}_1={}^{II}\!E^{(02)}_1=\nul$.
It follows that 
$^{I}\!E^{(20)}_2\cong{} ^{II}\!E^{(02)}_2$. We describe these $E_2$-terms in
a definition.

\definition{2.11 Definition} Let $\g M(\g A,\g B,P,Q,\pra~{},\prb~~)$ be Morita
context of \si\ Banach algebras, and let $D\colon\g A\to \g A^*$ be a
derivation. Then $D\in \Cal Z^1_{\g M}(\g A,\g A^*)$, if the tri-linear form
$(p,q,a)\mapsto\l\pra pq,D(a)\r$ is expressed as
$$
\l\pra pq,D(a)\r=\psi(ap,q)-\psi(p,qa),\quad p\in P,q\in Q, a\in \g A
$$
for some bilinear form $\psi\colon P\times Q\to \Bbb C$. We define a
subgroup of $\Ho1{\g A}$ by
$$
\HoM{\g A}=\Cal Z^1_{\g M}(\g A,\g A^*)\bigl/\Cal B^1(\g A,\g
A^*)\bigr.,
$$
where $\Cal B^1(\g A,\g A^*)$ is the space of inner derivations. The
subgroup $\HoM{\g B}$ of $\Ho1{\g B}$ is defined analogously.
\enddefinition

\example{2.12 Examples} We give some instances, where
$\HoM{\cdot}=\Ho1{\cdot}$
\flp1 Obviously (but importantly), if $\g A$ is \wae, then  $\HoM{\g
A}=\Ho1{\g A}=\nul$
\flp2 Suppose that $\g B$ is unital and that $\prb~~$ implements
an epimorphism $Q\pot{\g A} P\mapsto\g B$. Then we may write $\bold1_{\g B}=\sum
\prb{q_n}{q_n}$ and any derivation $D$ has $\l\pra
pq,D(a)\r=\psi(ap,q)-\psi(p,qa)$ with
$$
\psi(p,q)=\sum \l\pra{p_n}q,D(\pra p{q_n}\r.
$$
\flp3 Let $\g M(\g A,\g B,P,Q,\dots)$ be a Morita context and let
$\Cal M=\2x2{\g A}PQ{\g B}$. Let $\g N(\Cal M,\g A,\pmatrix
\g A\\Q\endpmatrix,(\g A\;P),\dots)$ be the iterated Morita context with pairings
defined by matrix multiplication. Then for any derivation $D\colon\Cal
M\to \Cal M^*$
$$
\multline
\l\pmatrix a\\q\endpmatrix(c\;p),D(\2x2 {a'}{p'}{q'}{b'}\r=\\
\l\2x2 a0q0,D(\2x2 cp00\2x2 {a'}{p'}{q'}{b'})\r-
\l\2x2 cp00\2x2 {a'}{p'}{q'}{b'},D(\2x2 a0q0)\r
\endmultline
$$
so that $\Cal H^1_{\g N}(\Cal M,\Cal M^*)=\Ho1{\Cal M}$
\endexample

Morita invariance of 1st degree cohomology now reads

\proclaim{2. 13 Theorem} Let $\g M(\g A,\g B, P, Q,\dots)$ be a full Morita
context of \si\ Banach algebras. Then
$$
\HoM{\g A}\cong\HoM{\g B}.
$$
\endproclaim
\medskip

\subheading{3. Applications to Banach algebras of approximable operators}

In [G1] we showed that for Banach algebras of approximable operators on
Banach spaces, whose dual have the bounded 
approximation property, Morita equivalence is determined by
approximate factorization. Specifically, $\A X$ and $\A Y$ are Morita
equivalent if and only if the multiplications $\Aa YX\pott\Aa XY\to\A
X$ and $\Aa XY\pott\Aa YX\to\A Y$ both are surjective. This can now be
improved to the much wider class for which $\A X$ is \si. As in [G1] the
key is to find the irreducible modules.
\proclaim{3.1 Lemma} $X$ is the only irreducible left module over $\A
X$.
\endproclaim
\demo{Proof} Let $L$ be a closed left ideal of $\A X$ and define
$$
\Psi(L)=\bigcup_{A\in L}A^*(X^*)\;.
$$
Then $\Psi(L)$ is a closed subspace of $X^*$: Let $\xi,\eta\in X^*$ and
$A,B\in L$. Choose $<x,x^*>=1$. Then
$$
\bigl((x\otimes\xi)A+(x\otimes\eta)B\bigr)^*x^*=A^*\xi+B^*\eta\;,
$$
so $\Psi(L)+\Psi(L)\subseteq\Psi(L)$. Clearly $\Bbb
C\Psi(L)\subseteq\Psi(L)$. To see that $\Psi(L)$ is closed, let
$A_n\in L,\; \xi_n\in X^*$ and suppose
$A_n^*\xi_n\to\xi$. Then $(x\otimes\xi_n)A_n\to x\otimes\xi\in L$, so
$\xi=(x\otimes\xi)^*x^*\in \Psi(L)$, i.e. $\Psi(L)$ is closed.

Now suppose that $\Psi(L)=X^*$ and let $x\otimes\eta$ be an arbitrary
rank-1 operator. Since $\Psi(L)=X^*$, we may choose $A\in L$ and
$x^*\in X^*$ so that $A^*x^*=\eta$. Then $x\otimes\eta=(x\otimes
x^*)A\in L$. Consequently $L$ contains all finite rank operators,
that is, $L=\A X$.

It is clear that $X$ is an irreducible $\A X$-module. Let $L$ be a
maximal modular left ideal of $\A X$. By the preceeding paragraph we
may choose $\xi\in X^*\setminus\Psi(L)$. Let $Q\colon\A X\to\A X/L$ be
the canonical map. Then
$$
x\mapsto Q(x\otimes\xi)\colon X\to\A X/L
$$
is a non-zero (bounded) module map between irreducible modules and
hence an isomorphism.
\enddemo

\proclaim{3.2 Theorem} Suppose that $\A X$ and $\A Y$ are
\si . Then $\A X$ and $\A Y$ are Morita equivalent, if and
only if multiplication implements surjections 
$$
\align
&\Aa XY\pott\Aa YX\to\A Y\\
&\Aa YX\pott\Aa XY\to\A X.
\endalign
$$
\endproclaim

\demo{Proof}  The proof is along the same lines as in [G1, Theorem 7.5] so we give
only a sketch. For brevity we set $\g A=\A X$ and $\g B=\A
Y$. Let the Morita equivalence be implemented by a full Morita context
$\g M(\g A,\g B, P, Q,\dots)$, and let $\Phi\colon P\pot{\g B}Q\to \g A$ be the
corresponding isomorphism. Since irreducibility of modules is Morita
invariant there is a bounded $\g B$-module isomorphism $\kappa\colon
Q\pot{\g A}X\to Y$. For each $q\in Q$ we define an operator $X\to Y$
by $x\mapsto \kappa(q\ot{\g A}x)$, and for each $p\in P$ we define
an operator $Y\to X$ by
$$y\mapsto \mu_X\circ(\Phi\ot{\g A}\bold
1_X)\circ(\bold 1_P\ot{\g B}\kappa^{-1})(p\ot{\g B}y).
$$
It follows that each $S\in \A X$ has a factorization
$$
S=\sum_nU_nV_n,\;U_n\in\Bb YX,V_n\in \Bb XY.
$$
Since $\A X$ is
\si\ we can actually obtain $U_n\in\Aa YX,V_n\in \Aa XY$.

Conversely, suppose that multiplications are surjective. Then $P=\Aa
YX$ and $Q=\Aa XY$ with pairings given by composition of operators
clearly satisfy the conditions of Lemma , so that the Morita context
$\g M(\g A,\g B, P, Q,\dots)$ is full.
\enddemo

We shall now investigate Hochschild cohomology groups related to
certain tensor products, of which the vector valued $L_p$-spaces are
prototypical. First a definition to describe the situation

\definition{3.3 Definition} Let $X$ and $Y$ be two Banach spaces.  We
say that {\it $Y$ is strongly finitely represented in $X$} if there is
$C>0$ such that for each finite dimensional subspace $S\subseteq Y$
there are a finite dimensional subspace $S\subseteq E$ and linear maps
$T_E \colon E\to X,\;\iota_E\colon X\to Y$ such that $\iota_ET_E$ is
the inclusion $E@>>>Y$ and $\nrm{T_E}\nrm{\iota_E}\leq C$.
\enddefinition

\example{3.4 Examples} The $\Cal L_p$-spaces are defined to be the Banach
spaces which are finitely represented in $\ell_p$ in the usual
(weaker) sense of finite representability. However they satisfy the
stronger condition of Definition 3.3, see [L\&R, Theorem III]. Any
Banach space is strongly finitely represented in the $C_p$-spaces of
W.B\. Johnson.
\endexample 

We now define the tensor products that are compliant with finite
representability. 

\definition{3.5 Definition} Let $Y$ and $Z$ be Banach spaces, let
$\nrm\cdot$ be a reasonable cross-norm on $Y\otimes Z$ with completion
$Y\btens Z$. For any complemented subspace $X\subseteq Y$, set $X\btens
Z=\cl(X\otimes Z)$. We say that $Y\btens Z$ is {\it compliable} if for
each complemented subspace $X\subseteq Y$ the map $x\otimes z\mapsto
Tx\otimes z$ defines a linear map $\Bb XY \mapsto \Bb{X\btens
Z}{Y\btens Z}$ of bound $\leq1$.
\enddefinition

\example{3.6 Examples} (1) Let $X$ be a complemented subspace of $Y$ with
complementation constant $C\geq1$ and
let $Y\btens Z$ be a compliable tensor product. Then $X\btens Z$ is
compliable tensor product and a
complemented subspace of $Y\btens Z$ with complementation constant $C$.  
\flp2 Let $X$ be any Banach space, and let $(\nu,
\Omega,\Sigma)$ be a measure space. Then the Banach spaces of
$p$-Bocher integrable $X$-valued functions $L_p(\nu,X),\;1\leq p\leq\infty$
are compliable tensor products.
\endexample

\proclaim{3.7 Proposition} Suppose that $Y$ is strongly finitely
represented in $X$ and that $Y\cong Y\oplus X$. Let $Z$ be any Banach
space and let $Y\btens 
Z$ be a compliable tensor product. If $\A{Y\btens Z}$ is \si, then
\roster
\item $\A{X\btens Z}$ is \si;
\item $\A{Y\btens Z}\sim\A{X\btens Z}$
\item $\Ho1{\A{Y\btens Z}}\cong\HoM{\A{X\btens Z}}$, where $\g M$
is the ca\-nonical Morita context $\g M(\A{Y\btens Z\oplus
X\btens Z},\A{X\btens Z},\dots)$.
\endroster
\endproclaim
\demo{Proof} Let for brevity $\g Y=Y\btens Z$ and $\g X=X\btens
Z$. Then $\g Y\cong\g Y\oplus\g X$ by compliance, so that $\A{\g X}$
may be viewed as a corner of $\A{\g 
Y}$. Since we have assumed that $\A{\g Y}$ is \si, any operator
$A\in\A{\g X}$ has an approximate factorization $A=\sum U_nV_n$ with
$U_n\in\Ff{\g Y}{\g X},\;V_n\in \Ff{\g X}{\g Y}$. By the defining
properties of strongly finite representability and of compliance, for
each $n\in\Bbb N$ there is a finite dimensional 
subspace $E_n\subseteq Y$ so that
$(\iota_{E_n}\btens\bold1_Z)(T_{E_n}\btens\bold1_Z)V_n=V_n$. It
follows that the multiplication $\A{\g X}\pot{\A{\g X}}\A{\g
X}\to\A{\g X}$ is surjective.

To show that it is injective, let $\phi$ be a bounded $\A{\g
X}$-balanced bilinear functional on $\A{\g X}$. According to the
decomposition $\g Y\cong \g Y\oplus 
\g X$ we may consider operators in $\A{\g Y}$ as matrices
$$
\2x2{\a11}{\a12}{\a21}{\a22},\quad \a11\in\A{\g Y},\a12\in\Aa{\g X}{\g
Y},\a21\in\Aa{\g Y}{\g X},\a22\in\A{\g X}. 
$$
Let $(\a ij),(\be ij)\in\F{\g Y}$ and define
$$
\multline
\widetilde\phi((\a ij),(\be ij))=\\
\lim_{E\to\infty}\phi(
(T_E\btens\bold1_Z\;\bold1_{\g X})(\a ij)\pmatrix\iota_E\btens\bold1_Z\\\bold1_{\g X}\endpmatrix,
(T_E\btens\bold1_Z\;\bold1_{\g X})(\be ij)\pmatrix\iota_E\btens\bold1_Z\\\bold1_{\g X}\endpmatrix),
\endmultline
$$
where the ordering is by inclusion of $E$'s (and formally $(T_E\btens\bold1_Z)\a
ij=0$ if the range of $\a ij$ is not included in $E\btens Z$.

This is easily seen to define a bounded $\A{\g Y}$-balanced bilinear
functional on $\A{\g Y}$ extending $\phi$ from the corner. Now suppose
that $\sum U_nV_n=0\;, U_n,V_n\in\A{\g X}$. Then $\sum\2x2
000{U_n}\2x2000{V_n}=0$. Since $\A{\g Y}$ is \si, we have that $\sum\widetilde\phi(\2x2
000{U_n},\2x2000{V_n})=0$. But then $\sum\phi(U_n,V_n)=0$,
i.e\. multiplication $\A{\g X}\pot{\A{\g X}}\A{\g
X}\to\A{\g X}$ is injective. 

To show (2) one may argue similarly or just appeal to Theorem 2.10(1). 

The statement (3) follows from Theorem 2.13 by noting that if $\g M$ is the
canonical Morita context of $\A{\g Y\oplus\g X}\sim\A{\g X}$, then $\HoM{\A{\g
Y\oplus\g X}}=\Ho1{\A{\g Y\oplus\g X}}$, so that $\Ho1{\A{\g
Y}}\cong\Ho1{\A{\g Y\oplus\g X}}=\HoM{\A{\g Y\oplus\g
X}}\cong\HoM{\A{\g X}}$.
\enddemo

In [G3] we proved that $\A{L_p(\nu,X)}$ is \wae\ for any infinite-dimen\-sional
$L_p(\nu)$ if $X$ has the bounded approximation property. This can now be improved to 

\proclaim{3.8 Corollary} Let $Y$ be an infinite-dimensional $\Cal
L_p$-space, $1\leq p<\infty,$ and let $\g Y=Y\btens Z$ be a compliable
tensor product. Then $\A{\g Y}$ is \wae, if and only if $\A{\g Y}$ is
\si. In particular for an infinite-dimensional $L_p(\nu)$,
$\A{L_p(\nu,X)}$ is \wae, if and only if it is \si.
\endproclaim

\demo{Proof} By [G3, Theorem 2.9] being \si\ is necessary. The rest follows
from noting that 
\flushpar - when $Y$is an  $\Cal L_p$-space, then $Y\cong
Y\oplus\ell_p$ ([L\&R, Theorem I]);
\flushpar - $L_p(\nu,X)$ is (isometrically) isomorphic to a compliable
tensor product $L_p(\nu)\btens X$; 
\flushpar - $\A{\ell_p(X)}$ is \wae, if and only it is \si, ([G3,
Corollary 4.2]). 
\enddemo

\medskip
We finish the paper with some clarification concerning direct sums and
\si ness. Recall that a
derivation $D\colon\A X\to\A X^*$ has the form $\l
A,D(B)\r=\tr((AB-BA)^*T)$ for an appropriate $T\in\Bn{X^*}$ and that
$D$ is inner if and only if $T\in\I{X^*}+\Bbb C\bold1_{X^*}$.

\medskip{\it Notation:\/} For brevity we shall in the following use
$[\;\pmb;\;]$ for commutators (of matrices).

\proclaim{3.9 Lemma} Let $D\colon\A{X\oplus Y}\to\A{X\oplus Y}^*$ be a
derivation implemented by 
$$
T=\pmatrix T_{11}&T_{12}\\
           T_{21}&T_{22}
  \endpmatrix
\in\Bn{X^*\oplus Y^*}.
$$
If $\A X\pot{\A{X}}\Aa YX=\Aa YX$, then $T_{12}\in\Ii{Y^*}{X^*}$. If $\Aa
XY\pot{\A{X}}\A X=\Aa XY$, then $T_{21}\in\Ii{X^*}{Y^*}$. If $\A{X\oplus
Y}$ is \si , then $T_{12}\in\Ii{Y^*}{X^*}$ and
$T_{21}\in\Ii{X^*}{Y^*}$. 
\endproclaim
\demo{Proof} For $A\in\F X$ and $P\in\Ff YX$ we have
$$
\align
\Vert A\Vert\Vert P\Vert\Vert D\Vert&\geq 
\vert<
\pmatrix A&0\\0&0\endpmatrix,
D(\pmatrix 0&P\\0&0\endpmatrix)>\vert\\
&=\vert \operatorname{tr}\biggl[\pmatrix A&0\\0&0\endpmatrix\pmb;\pmatrix
0&P\\0&0\endpmatrix\biggr]^*\pmatrix T_{11}&T_{12}\\T_{21}&T_{22}
\endpmatrix)\vert\\
&=\vert\operatorname{tr} \pmatrix
0&0\\0&(AP)^*T_{12}\endpmatrix\vert\\ 
&=\vert\operatorname{tr} ((AP)^*T_{12})\vert\;. 
\endalign
$$
Thus $(A,P)\mapsto\operatorname{tr} ((AP)^* T_{12})$ defines a
bounded, balanced 
bilinear form. Invoking the hypothesis we see that $T_{12}$ defines a
bounded linear form on $\Aa YX$, i.e\. $T_{12}$ is integral. The proof
of the second statement follows analogously. 
Now suppose that $\A{X\oplus Y}$ is \si .
This means that the four multiplications
$$
\align
&E_{11}\colon\pmatrix \A X&\Aa YX\\0&0\endpmatrix\pot{\A{X\oplus Y}}
\pmatrix \A X&0\\\Aa XY&0\endpmatrix\to\A X\\
&E_{12}\colon\pmatrix \A X&\Aa YX\\0&0\endpmatrix\pot{\A{X\oplus Y}}
\pmatrix 0&\Aa YX\\0&\A Y\endpmatrix\to\Aa YX\\
&E_{21}\colon\pmatrix 0&0\\\Aa XY&\A Y\endpmatrix\pot{\A{X\oplus Y}}
\pmatrix \A X&0\\\Aa XY&0\endpmatrix\to\Aa XY\\
&E_{22}\colon\pmatrix 0&0\\\Aa XY&\A Y\endpmatrix\pot{\A{X\oplus Y}}
\pmatrix 0&\Aa YX\\0&\A Y\endpmatrix\to\A Y
\endalign
$$
are linear topological isomorphisms. 
For $A\in\F X,\;P_1,P_2\in \Ff YX,\;B\in\F Y$ we get
$$
\align
(\Vert A\Vert+\Vert P_1\Vert)(\Vert P_2\Vert+\Vert B\Vert)\Vert
D\Vert&\geq\vert<\pmatrix A&P_1\\0&0\endpmatrix,D\pmatrix
0&P_2\\0&B\endpmatrix>\vert\\
&=\vert\operatorname{tr}(\biggl[\pmatrix A&P_1\\0&0\endpmatrix\pmb;\pmatrix
0&P_2\\0&B\endpmatrix\biggr]^*\pmatrix
T_{11}&T_{12}\\T_{21}&T_{22}\endpmatrix)\vert\\ 
&=\vert\operatorname{tr}(\2x2 0 0 {(AP_2+P_1B)^*}0\2x2
{T_{11}}{T_{12}}{T_{21}}{T_{22}})\vert\\
&=\vert\operatorname{tr}((AP_2+P_1B)^*T_{12})\vert 
\endalign
$$
Invoking the topological isomorphism $E_{12}$ we see that $T_{12}$
implements a bounded linear functional on $\Aa YX$ i.e\. $T_{12}$ is
integral. Similarly, $T_{21}$ is integral.
\enddemo

\proclaim{3.10 Theorem} Let $X$ and $Y$ be infinite dimensional Banach
spaces and assume that $\A X$ and $\A Y$ are WA. Then 
$\A{X\oplus Y}$ is WA if and only if $\A{X\oplus Y}$ is
\si\ and $\sup\{|\tr PQ|\mid P\in\Ff YX_1,\;Q\in \Ff XY_1\}=+\infty$.
\endproclaim
\demo{Proof} Since $X\oplus Y$ is decomposable, \si ness is
necessary ([G3, Theorem 2.9]). The other necessary condition is [B1,
Proposition 4.1].

 Now assume that $\A{X\oplus Y}$ is \si\ and that the supremum
is infinite. Then, by the lemma, $T_{12}$ and $T_{21}$ are integral
and since $\A X$ is WA, there
is $\lambda\in\Bbb C$ so that $T_{11}+\lambda\bold 1_{X^*}\in \I
{X^*}$. By subtracting the inner derivation given by
$$
\2x2{T_{11}+\lambda\bold 1_{X^*}}{T_{12}}{T_{21}}0
$$
we see, after adjusting by a multiple of $\bold 1_{X^*\oplus Y^*}$ if
necessary, 
that in order to show that $\A{X\oplus Y}$ is WA it suffices to 
look at derivations given by $T\in\Bn{X^*\oplus Y^*}$ of the form
$$
T=\2x2 000{T_{22}}.
$$
Since $\A Y$ is WA, there is $\lambda\in\Bbb C$ so that
$T_{22}+\lambda\bold1_{Y^*}\in\Cal I(Y^*)$. If $\lambda\neq 0$, then 
$$
\2x2 000{\bold1_{Y^*}}
$$
defines a bounded derivation. But then there is $K>0$ so that 
$$
|\tr(PQ)|=\bigl|\tr\bigl[\2x2 00Q0\pmb;\2x2 0P00\bigr]^*\2x2
000{\bold1_{Y^*}}\bigr|\leq K\Vert P\Vert\Vert Q\Vert
$$
contrary to assumption. Hence $T_{22}\in\Cal I(Y^*)$ and the
derivation is inner as wanted. 
\enddemo

\remark{Question} Theorem 6.10 of [G,J\&W] states that if $\A{X\oplus Y}$ is
amenable, then at least one of $\A X$ and $\A Y$ is amenable. Can \lq
amenable' be replaced by \lq\wae'? Note that $\A{X\oplus C_p}$ is \wae\
for any Banach space ([B1, corollary 3.3]), so that the conjecture is best possible.
\endremark

\enddocument